\documentclass[12pt]{amsart}
\input{epsf}
\usepackage{amssymb}
\usepackage{amsfonts}
\usepackage{amsmath}
\usepackage[matrix,arrow,graph]{xy}

\begin{document}
\title{From moment graphs to intersection cohomology}
\begin{abstract}
We describe a method of computing equivariant and ordinary
intersection cohomology of certain varieties with actions of algebraic 
tori, in terms of structure of the zero- and one-dimensional
orbits.  The class of varieties to which our formula applies 
includes Schubert varieties in flag varieties and affine flag 
varieties.  We also prove a monotonicity result on local intersection
cohomology stalks.
\end{abstract}

\author{Tom Braden}
\address{Harvard University}
\email{braden@math.harvard.edu}
\author{Robert MacPherson}
\address{Institute for Advanced Study}
\email{rdm@ias.edu}

\maketitle

\theoremstyle{plain}
\newtheorem{thm}{Theorem}[section]
\newtheorem{lemma}[thm]{Lemma}
\newtheorem{prop}[thm]{Proposition}
\newtheorem{cor}[thm]{Corollary}
\theoremstyle{definition}
\newtheorem*{conj}{Conjecture}
\newtheorem*{defn}{Definition}
\newtheorem*{nota}{Notation}
\newtheorem*{claim}{Claim}
\newtheorem*{rmks}{Remarks}

\hyphenation{Mac-Pher-son}
\hyphenation{canon-ical}

\newcommand{\ra}{\mathop{\rightarrow}\limits}
\newcommand{\la}{\mathop{\leftarrow}\limits}
\newcommand{\lra}{{\longrightarrow}}
\newcommand{\lla}{{\longleftarrow}}
\newcommand{\wt}{\widetilde}
\newcommand{\ol}{\overline}

\def\rank{\mathop{\rm rank}\nolimits}
\def\ker{\mathop{\rm ker}\nolimits}
\def\coker{\mathop{\rm Coker}\nolimits}
\def\im{\mathop{\rm Im}\nolimits}
\def\Id{{\rm Id}}
\def\id{{\rm id}}
\def\lcm{\mathop{\rm lcm}\limits}
\def\Hom{\mathop{\rm Hom}}
\def\Gr{\mathop{\rm Gr}\nolimits}
\def\Sym{\mathop{\rm Sym}}
\def\Aut{\mathop{\rm Aut}}

\newcommand{\C}{{\mathbb C}}
\newcommand{\Z}{{\mathbb Z}}
\newcommand{\R}{{\mathbb R}}
\newcommand{\Rge}{\R_{\ge 0}}
\newcommand{\PP}{{\mathbb P}}
\newcommand{\Q}{{\mathbb Q}}
\newcommand{\F}{{\mathbb F}}
\newcommand{\HH}{{\mathbb H}}
\newcommand{\cal}{\mathcal}
\newcommand{\cS}{{\cal S}}
\newcommand{\cV}{{\cal V}}
\newcommand{\cE}{{\cal E}}
\newcommand{\cF}{{\cal F}}
\newcommand{\cH}{{\cal H}}
\newcommand{\cL}{{\cal L}}
\newcommand{\cP}{{\cal P}}
\newcommand{\cA}{{\cal A}}
\newcommand{\cR}{{\cal R}}
\newcommand{\cD}{{\cal D}}
\newcommand{\cM}{{\cal M}}
\newcommand{\cN}{{\cal N}}
\newcommand{\IH}{{\cal I}{\cal H}}
\newcommand{\mft}{{\mathfrak t}}
\newcommand{\bdy}{{\partial}}

\newcommand{\IC}{{\mathbf {IC^{\textstyle \cdot}}}}
\newcommand{\mb}{\mathbf}

\section{Statement of the Main Results}

To a variety $X$ with an appropriate torus action (\S 1.1), we will 
associate a {\em moment graph} (\S 1.2), a combinatorial object which 
reflects the structure of the $0$ and $1$-dimensional orbits.  There is a 
canonical {\em sheaf} (\S 1.3) on the moment graph, combinatorially 
constructed from it (\S 1.4), which we denote by $\cM$.  
The main result (\S 1.5) uses the sheaf $\cM$ to compute
the local and global equivariant and ordinary 
intersection cohomology of $X$ functorially.

\subsection{Assumptions on the Variety $X$}\label{Assumptions on the variety}

We assume that $X$ is an irreducible complex algebraic variety 
endowed with two structures:

\begin{enumerate}
\item An action of an algebraic torus $T \cong (\C^*)^d$.  We assume that 
\begin{enumerate}\item 
for every fixed point $x\in X^T$ there is a one-dimensional subtorus
which is contracting near $x$, i.e.\ there is a homomorphism
 $i\colon \C^* \to T$ and a Zariski open neighborhood $U$ of $x$ so that 
 $\lim_{\alpha \to 0} i(\alpha)y = x$
for all $y\in U$ (this implies $X^T$ is finite), and
\item $X$ has finitely many one-dimensional orbits 
%(note that if $X$ is complete, then this  
%follows from (a)).
\end{enumerate} 
\item A $T$-invariant Whitney stratification by affine spaces. 
\end{enumerate}

It follows that each stratum contains exactly one fixed point,
since a contracting $\C^*$ action
on an affine space must act linearly with respect to 
some coordinate system
(see \cite{B-B}, Theorem 2.5).
Let $C_x$ denote the stratum containing the fixed point $x$, 
so $X = \bigcup_{x\in X^T} C_x$.  
Every one dimensional orbit $L$ has exactly two distinct 
limit points: the $T$ fixed point $x$ in the stratum $C_x$ containing 
$L$ and another fixed point lying in some stratum in the closure of $C_x$.  

The main case we are interested in is 
when $X$ is a Schubert variety in a flag variety or affine flag variety.
More generally, if $M$ is a smooth projective variety with a $T$ action
satisfying (a) and (b) above, one can take a homomorphism
$\C^* \to T$ for which $M^{\C^*} = M^T$ and consider the corresponding 
Bialnicki-Birula
decomposition of $M$ into cells.  If it is a stratification, then the closure 
of any cell satisfies our hypotheses.

\subsection{Moment graphs} \label{Moment graphs}
Let $\mft$ be a complex vector space.  A {\em $\mft$ moment graph} 
$\Gamma$ is a finite graph with a two additional structures:  
\begin{enumerate}
\item for each edge $L$, a one dimensional subspace $V_L$ of the dual vector space 
$\mft^*$ called the {\em direction} of $L$, and
\item a partial order $\leq$ on the set of vertices with 
the property that if an edge $L$ connects vertices $x$ and $y$, 
then either $x\leq y$ or $y\leq x$ (but $y\neq x$).
\end{enumerate}

We denote the set of vertices of $\Gamma$ by $\cV$, and the set of 
edges by $\cE$.  For a vertex $x \in \cV$, we denote by $U_x$ (for ``up'') 
the set of edges connecting $x$ to a vertex $y$ where $x\leq y$, and by 
$D_x$ (for ``down'') the set of edges connecting $x$ to a 
vertex $y$ where $y\leq x$.  

\vspace{.1in}
\noindent {\bf Constructing a moment graph from $X$}.  Given a variety 
$X$ as in \S 1.1, we construct a moment graph $\Gamma$ as follows.  
The vertices of $\Gamma$ are the $T$ fixed points in $X$, and 
the edges of $\Gamma$ are the one dimensional orbits of $X$.  
The vector space $\mft$ is the Lie algebra of $T$.  For an edge 
$L\in \cE$, every point on the one dimensional orbit has the same 
stabilizer in $T$; its Lie algebra is a hyperplane in $\mft$.  
The direction $V_L$ is the annihilator of that hyperplane in $\mft^*$.  
The partial order is defined by saying that for $x$ and $y$ in $\cV$, 
$x\leq y$ if and only if the stratum $C_y$ is in the closure of $C_x$.  
Note that $D_x\subset \cE$ is the set of one dimensional orbits 
contained in $C_x$.   

\begin{rmks}  Similar structures
(for smooth varieties) are considered by Guillemin and Zara in \cite{GZ}, 
\cite{GZ2}, \cite{GZ3}.

The term moment graph is motivated by the following.  If $X$ is projective,
%has a $T_K$-invariant symplectic form, where $T_K \cong (S^1)^d$ is the maximal
%compact subgroup of $T$.  
there is a moment map $\mu \colon 
X \to \mft^*_K$ to the dual of the Lie algebra of the maximal compact torus
$T_K \subset T$.  For $L\in \cE$, the image $\mu(\ol{L})$ is a
line segment joining $\mu(x)$ and $\mu(y)$, where $\{x, y\} = \ol{L}\cap \cV$.
  The vector $\mu(x) - \mu(y)$ 
spans the space $V_L$, using the identification 
$\mft^* \cong \mft^*_K \otimes_\R \C$.
\end{rmks}

\subsection{Sheaves on the moment graph} Let $A=\Sym(\mft^*)$ be 
the ring of polynomial functions on $\mft$.  Given $L\in \cE$, denote  
the quotient ring $A/V_L A$ by $A_L$.  For us, a ``module'' over $A$ or $A_L$ will
always be a finitely generated graded module. 

\begin{defn}  Let $\Gamma$ be a $\mft$ moment graph. 
A {\em $\Gamma$-sheaf} $\cM$ is a triple 
$\cM = (\{M_x\}, \{M_L\}, \{\rho_{x,L}\})$ where $M_x$ is an 
$A$-module defined for each vertex $x \in \cV$, $M_L$ is an $A_L$-module (also an 
$A$-module by the homomorphism $A \to A_L$) defined for 
each $L \in \cE$, and $\rho_{x,L}\colon M_x \to M_L$ is a 
homomorphism of $A$-modules defined whenever the vertex $x$ lies on the edge $L$.
\end{defn}

Let $S(\Gamma)$ be the finite set $S(\Gamma) = \cV \cup \cE$ 
of vertices and edges of $\Gamma$.  Given a subset $Z\subset S(\Gamma)$,
we define a module $\cM(Z)$ of ``sections'' on $Z$ by
\[\cM(Z) = \{(\{s_x\},\{s_L\}) \in \bigoplus_{a \in Z} M_a \mid
 \rho_{x,L}(s_x) = s_L \;\text{if $x$ lies on $L$}\;\}.\] 
In other words, an element of $\cM(Z)$ is a choice of an element of $M_x$ for each 
$x\in Z\cap \cV$, together with a choice of an element of $M_L$ for each 
$L\in Z\cap \cE$, such that these elements are compatible under the 
homomorphisms $\rho_{x,L}$.

In a similar way, we have a ``sheaf of rings'' $\cA = (\{A_x\},\{A_L\},\{q_{x,L}\})$
on $\Gamma$, given by letting $A_x = A$ for all $x \in \cV$, and letting
the $q_{x,L} \colon A_x \to A_L = A/V_LA$ be the quotient homomorphisms.  Then
we can define a ring of sections $\cA(Z)$ of $\cA$ in the same way as above, 
and 
$\cM(Z)$ becomes a module over $\cA(Z)$.

A $\Gamma$-sheaf $\cM$ can be thought of as a sheaf in the usual sense.  
Put a topology on $S(\Gamma)$ by 
declaring $O \subseteq S(\Gamma)$ to be open if whenever $x\in O\cap \cV$ is a vertex, 
all edges $L\in \cE$ adjacent to $x$ are in $O$ as well.
Given a $\Gamma$-sheaf $\cM$, sending an open set $O$ to 
$\cM(O)$ defines a sheaf on $S(\Gamma)$; restriction homomorphisms 
are defined in the obvious way.  
In the same way $\cA$ defines a sheaf of rings 
on $S(\Gamma)$, and the sheaf $\cM$ is a sheaf of modules over $\cA$.

\begin{prop}
This association gives a bijection between $\Gamma$-sheaves and sheaves of modules over 
$\cA$ on the topological space $S(\Gamma)$.
\end{prop}
Because of this, we will also refer to $\Gamma$-sheaves as $\cA$-modules.
\begin{proof}
If $\Sigma \subset S(\Gamma)$, we define $\Sigma^\circ$
to be the minimal open set with the same vertices as $\Sigma$.
If $\Sigma$ is a complete subgraph of $\Gamma$, then 
restriction gives an isomorphism
$\cM(\Sigma^\circ) = \cM(\Sigma)$.  

The proposition now follows immediately, since 
the $\Gamma$-sheaf can be recovered from the sheaf on 
$S(\Gamma)$ as follows: 
\[M_x = \cM(x^\circ),\;\, M_L =\cM(L),\] and $\rho_{x,L}$ is given 
by restriction $\cM(x^\circ) \to \cM(L)$.
\end{proof}
\subsection{Construction of the canonical 
$\Gamma$-sheaf $\cM$}\label{constructing M}

For an $A$-module $M$, we denote by  $\ol{M}$ the graded vector space 
$M \otimes_A \C = M/(\mft^*)M$.  Recall that a projective cover $P$ of 
an $A$-module $M$ is a free $A$-module on the smallest number of 
generators with a surjection $P \to M$.  This is equivalent to saying that 
the induced homomorphism 
$\overline{P} \to \overline{M}$ is an isomorphism.  

A projective cover $P$ may be constructed by setting $P=\ol{M} \otimes A$, 
and defining the map to $M$ by choosing any splitting 
of the quotient $M \to \ol{M}$.  In particular, projective covers of
$M$ are isomorphic up to a non-unique isomorphism.

Given a $\mft$ moment graph $\Gamma$ arising from a variety $X$, there is a 
canonical $\Gamma$-sheaf $\cM$ constructed by the following inductive 
procedure.  Begin at the ``top'' of 
$\Gamma$: since $X$ is irreducible by assumption, there is a unique vertex
$x_0$ which is maximal in the partial order; put $M_{x_0} = A$.  

Now suppose $\cM$ is known on the full subgraph
$\Gamma_{>x}$ consisting of all vertices $y$ with $y > x$, 
together with all edges joining them.  We extend it to 
$\Gamma_{\geq x}$, the full subgraph of all vertices $y\ge x$.  
First extend it to the set $\wt\Gamma_{>x} = \Gamma_{>x}\cup U_x = \Gamma_{\ge x}\setminus\{x\}$ 
as follows.  If $L\in U_x$ and $y \in \Gamma_{>x}$ is the upper vertex of $L$,
put $M_L = M_y/V_LM_y$ and let $\rho_{y,L}$ be the quotient homomorphism.

Next, define a module $M_{\bdy x}$ to be the image of the restriction 
homomorphism 
\begin{equation} 
\label{restriction}\phi\colon\cM(\wt\Gamma_{>x}) \to \cM(U_x).
\end{equation}
Then let $M_x$ be the projective cover of $M_{\bdy x}$.
The composition
\[M_x \to M_{\bdy x} \subset \cM(U_x) = \bigoplus_{L \in U_x} M_L\]
defines the homomorphisms $\rho_{x,L}$. 

Since projective covers are always isomorphic, this defines a sheaf
uniquely up to isomorphism.  To get a strong functorial result, 
we need to show our sheaves are ``rigid''.  This follows from the
following local result.

\begin{prop} \label{local noauto}
If $M_x \to M_{\bdy x}$ and $N_x \to M_{\bdy x}$ are two projective 
covers, then there is a unique isomorphism $M_x \to N_x$ commuting 
with the projective cover homomorphisms.
\end{prop}

The proof, which we postpone, uses the algebraic geometry of $X$.  
Denote by $\Aut_\cA(\cM)$ the group of automorphisms of $\cM$ as a graded 
$\cA$-module.

\begin{cor} \label{noauto}
The restriction $\Aut_\cA(\cM) \to \Aut_A(\cM(x_0))$ is a bijection, 
where the second group is the group of automorphisms of $\cM(x_0)$ as
a graded $A$-module.
As a result, the group of automorphisms of $\cM$ is just multiplication 
by scalars in $\C^*$.
\end{cor}

\vspace{.1in}
\noindent{\bf Another definition of $\cM$.}
There is another way to 
describe the sheaf $\cM$.  Call an $\cA$-module
$\cN$ {\em pure} if for all $x\in \cV$
\begin{enumerate}
\item $\cN(x)$ is a free $A$-module,
\item $\cN(L) = \cN(x)/V_L\cN(x)$ whenever $L \in D_x$, and
\item the restrictions of $\cN(x^\circ)\to \cN(U_x)$ and 
$\cN(\wt\Gamma_{> x}) \to \cN(U_x)$ have the same image.
\end{enumerate}

\begin{thm} $\cM$ is the unique 
indecomposable pure sheaf with $\cM(x_0) = A$.
Any pure $\cA$-module is a direct sum of sheaves     
obtained by applying the same construction to the subgraphs $\Gamma_{\le x}$
consisting of all vertices $y \le x$ and all edges joining them, or
by applying shifts to these sheaves.
\end{thm}

\begin{proof} The first statement is clear.  The second follows by an easy 
induction from the following statement: if $\cS$ is a pure $\cA$-module 
on $\Gamma$, and $\cS|_{\Gamma_{>x}}$ splits as a direct sum of two pure 
sheaves, then this splitting can be extended to $\Gamma_{\ge x}$ 
(see \cite{BBFK2}, Theorem 2.3
for the analogous result for toric varieties).

\end{proof}

\subsection{The main results}  \label{The main results} Suppose that a torus $T$ acts on a variety 
$X$ as in \S1.1, that the $\mft$ graph $\Gamma$ is constructed from 
$X$ as in \S1.2, and the $\Gamma$-sheaf $\cM$ is constructed from 
$\Gamma$ as in \S1.4.  
\begin{thm} There is a canonical identification 
\[IH_T^*(X)=\cM(\Gamma)\]
of the $T$-equivariant intersection cohomology of $X$ with the space of the 
global sections of $\cM$.  They are free $A$-modules.
The intersection cohomology of $X$ is given by
\[IH^*(X)=\ol{\cM(\Gamma)}=\cM(\Gamma)\otimes_A \C.\]
\end{thm}

The local intersection homology groups of $X$ at $x\in X$ are invariants of 
the singularity type of $X$ at $x$.  Since these are constant along a 
stratum $C_x\subset X$, to know them all it is enough to compute them at 
the fixed point $x\in C_x$.

\begin{thm}\label{local equivariant IH}
The local equivariant intersection cohomology at $x\in X$ is  
canonically isomorphic to the stalk $M_x$:
\[IH_T^*(X)_x=\cM(\{x\})=M_x\]
The local intersection cohomology of $X$ is given by
\[IH^*(X)_x=\ol{\cM(\{x\})}=\ol{M_x}.\]
\end{thm}

It follows from results in \cite{GKM} that similar calculations hold in 
ordinary cohomology if the sheaf $\cM$ is replaced by the sheaf $\cA$.  
We have $H_T^*(X)=\cA(\Gamma)$; $H^*(X)=\ol{\cA(\Gamma)}$; and 
(trivially) $H_T^*(X)_x=\cA(\{x\})=A$, and $H^*(X)_x=\ol{\cA(\{x\})}=\C$.

\begin{thm}
The module structures over the cohomology ring of the intersection 
cohomology groups mentioned above are given by the module structure 
over $\cA$ of $\cM$.  For example, the module structure of $IH^*(X)$ 
over $H^*(X)$ is the module structure of $\ol{\cM(\Gamma)}$ over 
$\ol{\cA(\Gamma)}$.
\end{thm}

Finally, we also prove an unrelated result, 
Theorem \ref{stalk monotonicity}, that says the intersection
cohomology stalks of $X$ can only grow larger at
smaller strata.  In the case of Schubert varieties, this 
gives another proof of an inequality on Kazhdan-Lusztig
polynomials originally proved by Irving \cite{I}.

\subsection{Remarks on the proof} \label{proof remarks} 
There is an equivariant intersection 
homology $\Gamma$-sheaf $\cM$ defined by  
\[M_x = IH^*_T(X)_x,\;\,M_L = IH^*_T(X)_L;\]
these are free modules over $A$, $A_L$ respectively.  
The homomorphism $\rho_{x,L}:M_x \to M_L$ is the composition
\[IH^*_T(X)_x\stackrel{\sim}{\leftarrow}IH^*_T(X)_{x\cup L}\to IH^*_T(X)_L.\]
We will prove the following slight improvement of Theorem 
\ref{local equivariant IH}:

\begin{thm}\label{mainthm}
The equivariant intersection homology $\Gamma$-sheaf is canonically isomorphic 
to the $\Gamma$-sheaf constructed in \S1.4.
\end{thm}

Using results of \cite{GKM}, this result implies all of the others in \S1.5
(the action of $T$ on $X$ is equivariantly formal, \cite{GKM}, for 
weight reasons).  Note also that because of Corollary \ref{noauto}, the 
identifications in section \S1.5 are all canonical.  Because of this,
we can use these sheaves to study how the intersection homology
sheaves extend each other to form more complicated perverse sheaves --
this will be explored in \cite{Bra}.

For the equivariant intersection homology $\Gamma$-sheaf, we 
have $M_{x_0}=A$ for the maximal vertex $x_0$ 
because $x_0$ is a smooth point of $X$.  
If $L \in D_y$, we have  $M_L = M_y/V_LM_y$ because 
$L$ and $y$ lie in the same stratum $C_y$.  So everything comes 
down to the calculation of $M_x$ in terms of the sheaf $\cM|_{\Gamma_{>x}}$.  

Let $N\subset X$ be a $T$-invariant normal slice to $C_x$ through $x$.
It can be embedded as a $T$-invariant closed subvariety of an affine space
$\C^r$ with a linear action of $T$.  If $T$ has a subtorus contracting
$\C^r$ to $\{0\}$, 
then $IH^*_T(N)$ is the projective cover of $IH^*_T(N_0)$ where 
$N_0 = N\setminus\{0\}$. This was originally proved by Bernstein and 
Lunts in \cite{BeL}; we will give a somewhat simpler proof.

Thus we see that our theorem amounts to showing that 
$M_{\partial x} = IH^*_T(N_0)$.
The localization theorem of \cite{GKM} says that for nice enough
(e.g. projective) $T$-varieties $Y$ with isolated fixed points,
restriction gives an injection $IH^*_T(Y) \subset IH^*_T(Y^T)$, and the
submodule is cut out by relations determined by the one-dimensional
orbits.  We can apply this to the projective variety $N_0/\C^*$ for a 
contracting subtorus $\C^* \subset T$; the result is that the
restriction  
\begin{equation} \label{xyz} 
IH^*_T(N_0) \to \bigoplus_{L\in U_x} IH^*_T(N_0)_L
\end{equation}
is an injection.  We then use a surjectivity result coming from 
the weight filtration of mixed Hodge
theory to argue that $M_{\bdy x}$ is the image of (\ref{xyz}).

This calculation is similar to the calculation of equivariant
IH for toric varieties described in \cite{BBFK}, \cite{BeL}, 
\cite{BrL2}.  In both cases there is an induction from larger strata
to smaller ones, at each step calculating the equivariant IH of
a punctured neighborhood of the singularity at a new stratum and taking
the projective cover.  

There are two main differences between our situation and the toric 
case. First, in our case we only need data
from the zero and one-dimensional orbits --- 
since the strata are contractible, knowledge of the stalk at 
$x$ is as good as knowledge of the stalks on all of $C_x$.
Second, in the toric case strata have affine neighborhoods
which are themselves unions of strata.  So the definition of
sheaves on fans, which are parallel to our $\Gamma$-sheaves, 
uses only one module for each stratum, and the computation of the
module analogous to $M_{\bdy x}$ is simpler.

We remark that the definition of $\cM$ in \S1.4 makes sense for general
moment graphs, whether or not they arise from a variety $X$.  However, 
we do not know how to show that it satisfies the degree restrictions of intersection 
cohomology, or that Corollary \ref{noauto} holds, without using the 
variety.  Such a result might be
useful, for instance, in interpreting the the coefficients of Kazhdan-Lusztig
polynomials for non-crystallographic Coxeter groups such as $H_3$ and $H_4$. 
In this case there is a clear definition of a moment graph (see \S\ref{G mod B}), 
but no underlying variety.

\subsection{Computational simplifications}
The main difficulty in computing the sheaf $\cM$ is in taking the
image of the homomorphism 
$\phi$ from (\ref{restriction}).  Fortunately, there is a 
major simplification, which we give as Theorem \ref{2planes}.
Essentially it says that to check whether an element of
$\cM(U_x)$ is in the image of $\phi$
it is enough to check that it can be extended 
to give sections of $\cM$ on planar subgraphs of $\Gamma_{>x}$. 

Some of the relations cutting out the image of $\phi$ are
easy to describe.  Suppose $x < y$, and take a subspace  $V\subset \mft^*$.
If we have an increasing path $x = x_1 < x_2 < \dots < x_n = y$
with $x_i$ joined to $x_{i+1}$ by an edge $L_i$, we call it a {\em $V$-path} if
$V_{L_i} \subset V$ for all $i$.

For an $A$-module $M$, we put
\[M_V = M \otimes_A (A/VA) = M/VM.\]  If we have a $V$-path from $x$ to $y$
as above, then $(M_{x_{i+1}})_V \to (M_{L_i})_V$ are isomorphisms, 
so we can compose their inverses with $(M_{x_i})_V \to (M_{L_i})_V$ 
to get a homomorphism $(M_{x})_V \to (M_{y})_V$.  
\begin{prop} \label{map independence}
This homomorphism depends only on $x$, $y$, and $V$, and not on the path.
\end{prop}
In particular, taking 
$V = \mft^*$, we get a homomorphism $m_{y,x}\colon\ol{M_{x}} \to \ol{M_{y}}$.

Similarly, by composing all but the first homomorphism we can get a 
homomorphism $m^V_{y,L_1}\colon (M_{L_1})_V \to (M_y)_V$.  It is
independent of the path chosen.
\begin{cor} \label{polygon rels}
If $\{\alpha_L\}_{L\in U_x}$ is in $M_{\bdy x}$, then 
$m^V_{y,L}(\alpha_L)$ is independent of $L$ and $V$ for any $y$.
\end{cor}

If there are only finitely many two dimensional orbits in the 
punctured neighborhood $N_0$ of a fixed point $x$, the image 
of the map (\ref{restriction}) is 
exactly the set of $\{\alpha_L\}$ satisfying these relations.
This will hold if and only if for every three distinct edges 
$L_1$, $L_2$, $L_3$ in $U_x$, the total span of the $V_{L_i}$
is three-dimensional.

This happens, for instance, when $X$ is a Grassmannian, i.e.\
$X = G/P$ where $P$ is a maximal parabolic in a 
semisimple complex algebraic group $G$.  In general, however, 
$N_0$ may have infinitely many two-dimensional orbits, and 
there are additional 
relations beyond those imposed by the corollary above.  
We will see an example of these extra relations
in \S\ref{sl2 example}.

\subsection{Acknowledgments} 

We would like to thank
Gottfried Barthel, Jean-Paul Brasselet, Karl-Heinz Fieseler, Ludger Kaup, 
Mark Goresky, Victor Guillemin, and Catalin Zara for stimulating 
conversations.

\section{Schubert varieties}

Our main motivation for this work was the case of Schubert varieties.  A flag 
variety $M$ is stratified by Schubert cells $C_x$, whose closures 
$\ol{C_x}$ are the Schubert varieties.  Our results give a functorial 
calculation of $IH^*(\ol{C_x})_y$ for a $T$-fixed point $y\in \ol{C_x}$.  
The Poincar\'e polynomials of these groups are the Kazhdan-Lusztig 
polynomials $P_{x,y}$, which are important in representation theory.  

Our calculation uses only data (the moment graph) from the 
interval $[y,x]$ in the Bruhat order.  Brenti \cite{Bre} has given a formula 
for the Kazhdan-Lusztig polynomials using only data from this graph (whereas 
the original Kazhdan-Lusztig algorithm used the entire interval $[0,x]$). We 
have not been able to understand Brenti's formulas in terms of our 
construction.
 
\subsection{Schubert varieties for a complex algebraic group.}\label{G mod B}

Let $G$ be a semisimple complex algebraic group, $B$ a Borel subgroup, $P 
\supseteq B$ a parabolic subgroup, and $T\subseteq B$ a maximal torus.  Then 
$M=G/P$ is a flag variety.  The Schubert cells $C_x$ of $M$ are the orbits of 
$B$ on $M$.  Let $X=\ol{C_x}\subseteq M$ be a Schubert variety.   Then the 
action of $T$ on $X$ satisfies the assumptions of \S \ref{Assumptions on the 
variety}, taking as strata the Schubert cells in $X$.  

To calculate the local or global intersection homology of $X$ as in 
\S\ref{The main results} we need to determine the moment graph for $X$, as 
defined in \S \ref{Moment graphs}.  Let $W$ be the Weyl group of $G$, and 
$W_P$ the parabolic subgroup of $W$ corresponding to $P$ ($W_P$ is the Weyl 
group of the Levi of $P$).  Then $W$ acts on $\mft^*$, the dual of the Lie 
algebra of $T$.  Let $v\in \mft^*$ be a vector whose stabilizer is $W_P$.  Then 
the following sets are canonically equivalent, and we abuse notation by 
identifying them:  the orbit ${\cal O}$ of $v$ under $W$, the quotient set 
$W/W_P$, the set of Schubert cells of $M$, and the set of fixed points $M^T$ 
of $M$.  There is a Bruhat partial order on this set (given by the usual 
Bruhat order on the maximal elements of the cosets of $W/W_P$), 
which corresponds to the closure relation on the Schubert cells.  
The moment graph $\Gamma$ of $X$ is determined as follows:

\begin{itemize}
\item The vertices of $\Gamma$ are those $y\in {\cal O}$ such that $y\leq x$.  
\item Edges $L$ connect pairs of vertices $y$ and $z$ such that $y=Rz$ where 
$R$ is a reflection (not necessarily simple) in $W$.
\item The direction $V_L \subset \mft^*$ is spanned by $y-z$.
\item The partial order is the Bruhat order.
\end{itemize}

So the embedding of ${\cal O}$  in $\mft^*$ gives a linear map of the moment 
graph to $\mft^*$ in which the direction of $L$ is the angle of the image of 
$L$.  Such a graph is drawn below in \S \ref{sl2 example}.

\subsection{Affine Schubert varieties and the loop group.}

Let $G$ be a semisimple complex algebraic group, $G(\C((t))\,)$ the 
corresponding loop group, $I$ an Iwahori subgroup, $P \supseteq I$ a 
parahoric subgroup.  Then $M=G/P$ is an affine flag variety.  The Schubert 
cells $C_x$ of $M$ are the orbits of $I$ on $M$.   Let $X=\ol{C_x}\subseteq 
M$ be an affine Schubert variety.  It is a finite dimensional projective 
algebraic variety, even though $M$ is infinite dimensional. Let $A\subseteq 
G(\C)$ be a maximal torus whose inclusion in $G(\C((t))\,)$ lies in $I$.  Let 
$T$ be the torus $A\times \C^*$ which acts on $M$ as follows:  $A$ acts 
through $G(\C)$ and $\C^*$ acts by ``rotating the loop", i.e. $\lambda\in 
\C^*$ sends the variable $t$ to $\lambda t$.  Then $T$ preserves $X$, and the 
action of $T$ on $X$ satisfies the assumptions of \S \ref{Assumptions on the 
variety}, taking as strata the Schubert cells in $X$.  

As before, to calculate the local or global intersection homology of $X$ we 
need to specify the moment graph for $X$.  Let $W$ be the affine Weyl group 
$W$ of $G(\C((t))\,)$, and $W_P$ the parabolic subgroup of $W$ corresponding 
to $P$ (note that $W_P$ is a finite group).  Then $W$ acts on $\mft^*$, the dual of 
the Lie algebra of $T$ in a somewhat nonstandard way satisfying the following 
properties:  
\begin{enumerate}
\item  The projection of $\mft^*$ to ${\mathfrak a}^*$ is $W$ equivariant, 
where the action of $W$ on ${\mathfrak a}^*$, the dual to the Lie algebra of $A$, is the 
standard one.  
\item  Reflections in $W$ act by pseudoreflections on $\mft^*$, 
i.e.\ order two affine maps that fix a hyperplane.  
\end{enumerate}
Up to affine equivalence, 
there are only two actions satisfying these properties, and the action in 
question is the one that is not the product action.  

With this set-up, the construction of $\Gamma$ is identical to the 
construction for semisimple algebraic groups above.  Let $v\in \mft^*$ be a 
vector whose stabilizer is $W_P$.  We identify the following sets, which are 
canonically equivalent:  the orbit ${\cal O}$ of $v$ under $W$, the quotient 
set $W/W_P$, the set of Schubert cells of $M$, and the set of fixed points 
$M^T$ of $T$.  There is a Bruhat partial order on this set, defined as above, 
which corresponds to the closure relations of the Schubert cells.  The moment 
graph $\Gamma$ of $X$ is determined by the same procedure: The vertices of 
$\Gamma$ are those $y\in {\cal O}$ such that $y\leq x$;   edges $L$ connect 
pairs of vertices $y$ and $z$ such that $y=Rz$ where $R$ is a reflection in 
$W$; the direction $V_L \subset \mft^*$ is spanned by $y-z$; and the partial 
order is the Bruhat order.  

As before, the embedding of ${\cal O}\cap X$  in $\mft^*$ gives the 
structure.  The points of $\cal O$ lie on a paraboloid in $\mft^*$.  The case 
of the loop Grassmannian (an affine flag manifold for a particular parahoric 
$P$), is worked out in \cite{AP}, which also has some pictures of $\cal O$. 

\subsection{Example}\label{sl2 example}

\begin{figure}
\begin{center}
\leavevmode
\hbox{%
\epsfxsize = 1.9in
\epsffile{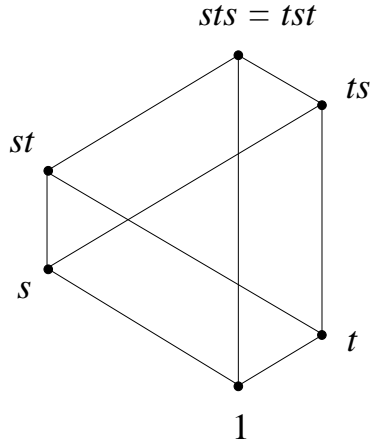}}
\end{center}
\caption{The moment graph for $G/B$, $G = SL_3$}
\label{hex}
\end{figure}

Take $G = SL_3(\C)$, and take 
$X = G/B$.  The moment graph is pictured in Figure \ref{hex}.
Since $X$ is smooth, we must have $M_w = A$ for all $w$.  Still, it is
instructive to see what Theorem \ref{mainthm} says in this case.

The induction begins with $M_{w_0} = A$ for the longest word $w_0 = sts$.
If $w = st$ or $ts$, there is only one edge $L$ in $U_w$, giving
$M_{\bdy w} = A/V_LA$. Since this module is generated in degree $0$, we have
$M_w \cong A$.   

If $w = s$ or $t$, there are two edges, say $L$ and $L'$, in $U_w$.
$M_{\bdy w}$ consists of pairs of polynomials in $M_L$ and $M_{L'}$
whose constant terms agree --- this is 
exactly the relation implied by Corollary \ref{polygon rels}.  
As a module this is just $A/V_LV_{L'}A$,
which again has a single generator in degree zero.

We see a new phenomenon when we look at $w = 1$.  The relations
from Corollary \ref{polygon rels} only affect the degree zero part; 
without further relations
we would have $\dim (M_{\bdy w})_2 = 3$, which would imply that $M_w$
has a generator in degree $2$.

The fact that we get the right relation from Theorem \ref{mainthm}
follows from the projective dual of Pappus' theorem in projective plane geometry.  
Degree two sections of $\cM = \cA$ can be seen as infinitesimal motions of the vertices
which preserve the parallelism classes of the edges.  If you 
remove the vertex labeled $1$ from Figure \ref{hex} and allow the remaining vertices
to move, the constraints imposed by the edges ensure that the three lines through $1$ 
will still meet in a point.

\subsection{}
For a Schubert variety $X \subset G/B$, there
is another description of $M_X= IH_T^*(X)$ as a module over $R=H^*_T(G/B)$,
coming from results due to Soergel (\cite{Soe2}, see \cite{Soe} for a 
non-equivariant version).  
In essence, he shows how to compute the equivariant cohomology of a 
resolution of $X$; by the decomposition theorem it is a direct sum 
of $M_X$ and shifted copies of $M_{X'}$ for smaller Schubert varieties $X'$.  
He proves that the $M_X$ are indecomposable $R$-modules, so in 
principle it is possible to compute the desired submodule.

\section{Equivariant intersection cohomology}
\subsection{Definitions and conventions} \label{IHT defs}
All our sheaves and 
cohomology groups will be taken with complex coefficients. For $X$ a complex
algebraic variety, let $\IC(X)$ be the intersection cohomology object in the
derived category $D^b(X)$, 
shifted so that it restricts to the constant local system in degree $0$
on the smooth locus of $X$. Its hypercohomology $\HH^d(\IC(X)) = IH^d(X)$ 
is the intersection cohomology of $X$.  

If $i\colon Y\to X$ is the inclusion of a subvariety, we put 
\[IH^d(X)_Y = \HH^d(i^*\IC(X)).\]  The adjunction 
$\IC(X) \to i_*i^*\IC(X)$ gives rise to a homomorphism
$IH^*(X) \to IH^*(X)_Y$.
If $i$ is a normally nonsingular inclusion, then there is a canonical
isomorphism $i^*\IC(X) \cong \IC(Y)$, giving a homomorphism
\[IH^*(X) \to IH^*(Y).\]

\subsection{Equivariant IH}
Now suppose an algebraic torus $T$ acts on $X$.  More sophisticated treatments
of equivariant intersection cohomology can be found in 
\cite{Bry},\cite{BeL},\cite{J}, but the following is enough for our purposes.
Fix an isomorphism $T \cong (\C^*)^d$, and let 
$E_k = (\C^k \setminus \{0\})^d$ carry the $T$-action given by termwise
multiplication. Let $E_k$ sit inside $E_{k+1}$ as the set of points
whose $(k+1)$st coordinates are all zero.  The quotient 
$B_k = E_k/T \cong (\C\PP^{k-1})^d$ is a finite-dimensional approximation to
the classifying space $BT = \bigcup B_k$.  The cohomology ring $H^*(BT)$ is
canonically isomorphic to the symmetric algebra $A = S(\mft^*)$.

Let $X_k = (X \times E_k)/T$.  The inclusion $X_k \subset X_{k+1}$ is
normally nonsingular. This gives a natural homomorphism
$IH^n(X_{k+1}) \to IH^n(X_k)$; it is an isomorphism when $2(k-1) \ge n$.  We
define the equivariant intersection cohomology by 
\[IH^n_T(X) = \lim_{\leftarrow} IH^n(X_k).\]
Similarly, if  $Y \subset X$ is a  $T$-invariant
subvariety, 
we put \[IH^n_T(X)_Y = \lim_{\leftarrow} IH^n(X_k)_{Y_k}.\] 

The natural map $X_k \to B_k$ makes
$IH^*(X_k)$ into a module over $H^*(B_k)$.  Taking limits, 
$IH^*_T(X)$ becomes an $A$-module.  Since $X_1 = X$, we have a 
homomorphism $IH_T^*(X) \to IH^*(X)$; it factors to give a
homomorphism \[\ol{IH^*_T(X)} \to IH^*(X).\] 

The following lemma gives the isomorphism used 
in \S1.6 to define the homomorphisms
$\rho_{x,L}$ in the sheaf $\cM$.
\begin{lemma} \label{contr}
Suppose $X$ has an algebraic $\C^*$ action which commutes with
the action of $T$ and contracts a locally closed subvariety $Y \subset X$ onto 
another subvariety $Y'$.  
Then $IH^*_T(X)_Y \to IH^*_T(X)_{Y'}$ is an isomorphism.
\end{lemma}

\subsection{Localization} We recall the result from \cite{GKM} that we
will need.

\begin{thm} \label{GKMloc} If either (a) $X$ is projective, or
(b) $IH^i(X)$ vanishes for $i$ odd,  then $IH^*_T(X)$ is a free
$A$-module, and the localization homomorphism 
\[\lambda\colon IH^*_T(X) \to IH^*_T(X)_{X^T}\] is an injection.

If $X$ has finitely many one-dimensional orbits, then 
\[(s_x) \in \bigoplus_{x\in X^T} IH^*_T(X)_x = IH^*_T(X)_{X^T}\]
is in the image of $\lambda$ if and only if $s_x$ and $s_y$ map
to the same element of $IH^*_T(X)_L$ whenever  
$L$ is a one-dimensional orbit with $\{x,y\} \subset \ol{L}$. 
\end{thm}

\subsection{Hodge Intersection Cohomology}
The proofs of our results rely on the weight filtration
on intersection cohomology of a complex variety.  This filtration
was constructed by Saito as part of his theory of mixed Hodge modules;
see \cite{S} for a good introduction.  In this section 
we extract some simple results from the theory which suffice for
our needs.

Given a complex variety $X$ and an open subvariety $U$, the
weight filtration is an
increasing filtration $W_iIH^*(X)$
and $W_iIH^*(X,U)$ on the intersection cohomology groups
$IH^*(X)$ and $IH^*(X,U)$.
It is strongly compatible with the homomorphisms in the long exact
sequence for the pair $(X,U)$: taking the associated graded
$\Gr^W_k$ of all terms in the sequence gives another exact
sequence.  

The homomorphisms $IH^*(X) \to IH^*(Y)$ induced by a normally nonsingular 
inclusion $X \to Y$ are also strictly compatible with the weight 
filtration.  Thus we get an induced weight filtration on the 
equivariant intersection cohomology $IH^*_T(X)$ of a variety with 
an algebraic $T$-action.

\begin{lemma} \label{semipurity}
We have $W_kIH^d(X,U) = 0$ if $k<d$.
\end{lemma}
\begin{proof} We have $IH^*(X, U) = \HH(i^!\IC(X))$, where 
$i\colon X\setminus U \to X$ is the inclusion.  According to 
\cite{S}, the functors $i^!$ and $\HH$ can only increase
weights.
\end{proof}

\begin{defn} We define the Hodge intersection cohomology of 
$(X,U)$ by 
\[HIH^d(X,U) = \Gr_d^WIH^d(X,U) = W_dIH^d(X,U).\]  
If $X$ carries a $T$-action, we let $HIH_T^d(X) = W_dIH_T^d(X)$.
\end{defn}

If $HIH^*(X) = IH^*(X)$, we say that $X$ is pure.
Projective varieties and quasiconical affine varieties are always pure.
The following proposition generalizes these two cases.

\begin{prop} \label{purity} 
If $X$ has an action of $T = \C^*$ which contracts $X$ onto
$X^T$, and $X^T$ is proper, then $X$ is pure. 
\end{prop}
\begin{proof} We have $IH^*(X)= \HH^*(\IC(X)) = \HH^*(i^*\IC(X))$, 
where $i\colon X^T \to X$ is the inclusion.  The middle expression
vanishes for weights less than the degree, while the right one vanishes
for weights greater than the degree (the functor $i^*$ can only decrease weights, and
hypercohomology of mixed Hodge modules on proper varieties preserves weights).
\end{proof}

\begin{thm} \label{Hodge monotonicity} If $U\subset X$ is an open 
subvariety, then the restriction $HIH^*(X) \to HIH^*(U)$ is a 
surjection.  If $X$
carries an action of $T$ and $U$ is $T$-invariant, then
$HIH^*_T(X)\to HIH^*_T(U)$ is a surjection.
\end{thm}

For example, take $X = \C\PP^1$, and $U = X \setminus \{p,q\}$ for $p\ne q$.
Then $IH^1(X) \to IH^1(U)$ is not surjective, but $HIH^1(U) = 0$.

\begin{proof} 
Lemma \ref{semipurity} implies that 
the connecting homomorphism $\Gr^W_dIH^d(U) \to \Gr^W_dIH^{d+1}(X,U)$ vanishes.
The equivariant case then follows from the nonequivariant case, since for
$k$ large enough we have $HIH_T^n(X) = HIH^n(X_k)$ and 
$HIH_T^n(U) = HIH^n(U_k)$.
\end{proof}

\subsection{Monotonicity for local stalks}
Theorem \ref{Hodge monotonicity} has the following consequence, which 
is independent from the rest of the paper.  Let $X$ be a $T$-variety 
satisfying the conditions of \S\ref{Assumptions on the variety}.  
Take $x,y \in X^T$ and assume that $x \le y$, so $x\in \ol{C_y}$.

\begin{thm}
\label{stalk monotonicity} There is a surjection $IH^*(X)_x\to IH^*(X)_y$.
\end{thm}
\begin{proof} For any $x\in X^T$, let $U$ be a $T$-invariant affine
neighborhood of $x$, and let $\rho_x$ be the 
composition of restriction and localization homomorphisms 
\[IH^*(X) \to IH^*(U) \stackrel{\sim}{\to} IH^*(U_x)_x = IH^*(X)_x.\]  It is
a surjection, using Theorem \ref{Hodge monotonicity} and Lemma \ref{contr}.
So it is enough to find a homomorphism $m\colon IH^*(X)_x\to IH^*(X)_y$ with
$m\rho_x = \rho_y$.
Such an $m$ is given by the composition 
\[IH^*(X)_x \stackrel{\sim}{\la} IH^*(U) \to IH^*(X)_{y'} 
\stackrel{\sim}{\la} IH^*(X)_{C_y} \stackrel{\sim}{\to} IH^*(X)_y,\] 
where $y'$ is any point in $C_y\cap U$. 
The last two isomorphisms result from the equisingularity of $X$ along $C_y$.
\end{proof}

The homomorphism
 $m$ does not depend on the choice of point $y'$, and in fact it
can be described in our moment graph language; it is
the homomorphism $m_{y,x}\colon \ol{M_x} \to \ol{M_y}$ given 
by Proposition \ref{map independence}.

If $X$ is a Schubert variety in a flag variety or affine flag variety, 
this gives an inequality on Kazhdan-Lusztig polynomials: let $P^i_{x,y}$ be the
$i$th coefficient of $P_{x,y}$.  
\begin{cor} $P^i_{x,z} \ge P^i_{y,z}$ if $x \le y$.
\end{cor}

This was proved algebraically in the case of ordinary flag varieties
by Irving (\cite{I}, Corollary 4), 
using the Koszul dual interpretation of
Kazhdan-Lusztig polynomials as multiplicities of simple objects in 
the socle filtration of a Verma module.  To our knowledge the corresponding
statement for affine flag varieties was not previously known.

\subsection{The local calculation}
The following theorem describes the local $IH^*_T$ groups of quasihomogeneous
singularities.  It was proved by Bernstein and Lunts in \cite{BeL}; we
will give a proof we feel is slightly simpler.
 
Suppose that a torus $T$ acts linearly on $\C^r$, and a subtorus
$\C^* \subset T$ contracts $\C^r$ to $\{0\}$.  Let $X\subset \C^r$
be a $T$-invariant variety, and put $X_0 = X\setminus\{0\}$. 
By Lemma \ref{contr}, we have $IH^*_T(X)_x \cong IH^*_T(X)$.
 
Recall that for any $A$-module
$M$ we put $\ol{M} = M \otimes_A \C$.
\begin{thm} \label{lociht} The restriction homomorphism makes
$IH^*_T(X)$ into a projective cover of $IH^*_T(X_0)$.
Its kernel is isomorphic to the local equivariant
intersection homology with compact supports
\[IH^*_{T,c}(X) = IH^*_T(X, X_0);\] it is a free $A$-module, and
$\ol{IH^*_{T,c}(X)} = IH^*_c(X)$. 
\end{thm}
\begin{proof} The claim that $IH^*_T(X)$ is a free module
follows from the collapsing of 
the spectral sequence $H^p(\bigcup B_k; IH^q(X)) \implies IH^{p+q}_T(X)$.  
This in turn happens because the intersection cohomology of 
the varieties $X_k = (X \times E_k)/T$ is pure, by Proposition \ref{purity}.
A similar argument shows that $IH^*_{T,c}(X)$ is free.

%Note that if $IH^*(X)$ vanishes in odd degrees, this spectral sequence
%collapses automatically.  This will be the case for all singularities

Theorem 9.1 of \cite{BeL} shows that 
$IH^*_T(X_0) = IH^*_{T/\C^*}(X_0/\C^*)$. 
Since $X_0/\C^*$ is projective, we have
$HIH^*_T(X_0) = IH^*_T(X_0)$, and so
Theorem \ref{Hodge monotonicity} implies that $IH^*_T(X) \to IH^*_T(X_0)$
is a surjection.

All that remains to prove the first statement 
is to show that $\ol{IH^*_T(X)} \to \ol{IH^*_T(X_0)}$ is
an injection.  But in the commutative square:
\[\xymatrix{ \ol{IH^*_T(X)} \ar[r]\ar[d] & \ol{IH^*_T(X_0)} \ar[d]\\
IH^*(X) \ar[r] & IH^*(X_0)
}\]
the left homomorphism 
is an isomorphism, and the lower one is an injection ---
it is an isomorphism in degrees $< \dim_\C(X)$ and $IH^*(X)$ vanishes
in higher degrees.

Finally, the second statement of the theorem follows from the 
long exact sequence for $IH^*_T$ of the pair $(X,X_0)$.
\end{proof}

\section{Proofs}
\subsection{The main theorem} 
\label{mtpf} We now have all the ingredients to prove  
Theorem \ref{mainthm}.  Let $\cM$ be the $\Gamma$-sheaf defined by the
inductive construction of 
\S\ref{constructing M}, and let $\cM'$ be the equivariant
intersection cohomology $\Gamma$-sheaf of \S\ref{proof remarks}.
We need to show that these sheaves are isomorphic.

Their restrictions to the top vertex $x_0$ clearly agree, since 
$M_{x_0} =A$ by definition and $X$ is smooth at $x_0$.  Further,
if $\cM$ and $\cM'$ agree at a vertex $y$, they agree at all 
edges $L\in D_y$, since $M_L = M_y/V_LM_y$ and $X$ is equisingular 
along $y \cup L \subset C_y$.

Now take a vertex $x \in \cV$ and assume inductively that $\cM$ and
$\cM'$ have isomorphic restrictions to $\Gamma_{> x}$.    
We need to show they agree on all of $\Gamma_{\ge x}$.  By the previous 
remark, they agree on $\wt\Gamma_{>x} = \Gamma_{>x} \cup U_x$. 

By our assumptions, $x$ has a $T$-invariant affine
neighborhood $U$.
\begin{lemma} There is a $T$-invariant closed subvariety $N \subset U$
which is a normal slice to $C_x$ through $x$.
\end{lemma} 
\begin{proof} We can find a diagonal linear action of $T$ on some
affine space $\C^r$, and an equivariant embedding 
$U\subset \C^r$. The tangent space $T_xC_x$ will be generated by
a subset of the coordinate directions.  Take the linear span of the
remaining coordinates and intersect with $U$.
\end{proof}

Applying Lemma \ref{contr} to $U$ and to $N$ and using 
the fact that $N\subset U$ is a normally nonsingular inclusion, 
we obtain isomorphisms $IH^*(U) \cong  IH^*(U)_x \cong 
IH^*(X)_x$ and $IH^*(U) \to IH^*(N)$.

If we can show that
$IH^*_T(N_0)\cong M_{\bdy x}$, where $N_0 = N\setminus \{x\}$, then 
Theorem \ref{lociht} implies that $M_x \cong IH^*_T(X)_x$, as required.  
Let $X_{>x} = \bigcup_{y>x} C_y$.  
Consider the following diagram of restriction
homomorphisms: 
%(where we use the isomorphisms 
%$IH^*_T(X)_N \cong IH^*_T(N)$, $IH^*_T(X)_{N_0} \cong IH^*_T(N_0)$):
\[\xymatrix{IH^*_T(X) \ar[r]\ar[d] & IH^*_T(X_{>x}) \ar[r]^{\alpha}
\ar[d]^{\gamma} &
 \cM(\wt\Gamma_{>x}) \ar[d]^{\phi}\\
  IH^*_T(N) \ar[r] & IH^*_T(N_0) \ar[r]^{\beta} & \cM(U_x)}\]
We will show that $\alpha$ is an isomorphism, $\beta$ is an injection, and
$\gamma$ is a surjection; the result follows.

To see that $\alpha$ is an isomorphism, we apply  
Theorem \ref{GKMloc} to $X_{>x}$; we claim that 
$IH^d(X_{>x})$ vanishes for $d$ odd.  It follows from the previous 
steps of the induction that the local stalks $IH^*(X)_y \cong \ol{M_y}$
vanish in odd degrees for $y > x$.  Applying Verdier duality 
we see that $\HH^*(i_y^!\IC(X))$ vanishes in odd degrees if
$i_y\colon C_y\to X$ is the inclusion and $y>x$.  Since this 
can also be written $IH^*(X_{>y}\cup C_y,X_{>y})$, the 
claim follows by induction using the long exact sequence of a pair.

To see that $\beta$ is an injection, note that any 
contracting subtorus $\C^*\subset T$ acts almost freely 
(only finite stabilizers) on $N_0$.  By \cite{BeL}, Theorem 9.1,
we have isomorphisms \[IH^*_T(N_0) \cong IH^*_{T/\C^*}(N_0/\C^*),\]
\[\cM(U_x) \cong \bigoplus_{y \in (N_0/\C^*)^T} IH^*_{T/\C^*}(N_0/\C^*)_y.\]
%where the sum on the right is over all fixed points of $N_0/\C^*$.  
Since $N_0/\C^*$ is projective, we can apply Theorem \ref{GKMloc}.

Finally, $\gamma$ is a surjection because 
$IH^*_T(X) \to IH^*_T(N)$ and $IH^*_T(N) \to IH^*_T(N_0)$
are surjections.  The first homomorphism factors as
$IH^*_T(X) \to IH^*_T(U)\stackrel{\cong}{\to} IH^*_T(N)$, so 
the surjectivity follows from Theorem \ref{Hodge monotonicity} and
Proposition \ref{purity}.  The second surjection is part of Theorem
\ref{lociht}. 

Note that we have shown that $IH^*(X)$ vanishes in odd degrees,
so Theorem \ref{GKMloc} can be applied to deduce the 
theorems in \S1.5 from Theorem \ref{mainthm}.

\subsection{Automorphisms} 
Proposition \ref{local noauto} now follows from 
Theorems \ref{mainthm} and \ref{lociht}, the
degree vanishing conditions for local intersection cohomology
and compactly supported intersection cohomology, and
the following lemma.

\begin{lemma} Let $M_i$, $M'_i$, $i= 1, 2$ be graded modules over a polynomial
ring $A$, with $M_i$ free, and let $\phi_i\colon M_i \to M'_i$ be 
homomorphisms with $\ol{\phi}\colon \ol{M_i} \to \ol{M'_i}$ an isomorphism.
Also suppose that for some $d\in \Z$, $M_1$ is generated in degrees $< d$
and $\ker \phi_2$ is generated in degrees $\ge d$.

Then if $f'\colon M'_1 \to M'_2$ is a homomorphism of graded modules, 
there is a unique 
$f\colon M_1 \to M_2$ so that $\phi_2 f = f'\phi_1$.
\end{lemma}

\subsection{Planar relations} \label{planar}
For the results of this last section, we need to assume that the moment
graph $\Gamma$ is constructed from a projective variety $X$.

Fix a vertex $x\in \cV$ of the moment graph.  If $H \subset \mft^*$
is a sub-vector space, consider the graph with the
same vertex set as $\Gamma$, but with only those 
edges $L$ of $\Gamma$ for which $V_L \subset H$.  Denote by $\Gamma^H$
the connected component of this graph containing $x$, and let
$\wt\Gamma^H_{>x}$, $U^H_x$ be the intersections of $\wt\Gamma_{>x}$, $U_x$
with $\Gamma^H$.

Let $\phi\colon\cM(\wt\Gamma_{>x}) \to \cM(U_x)$ and  $\phi^H\colon
\cM(\wt\Gamma^H_{>x}) \to \cM(U^H_x)$ be the restrictions.   Given
$\xi \in \cM(U_x)$, let $\xi^H$ be its restriction to $\cM(U^H_x)$.
Let $\cH$ be the set of all two-dimensional subspaces $H\subset
\mft^*$ for which $\Gamma_{>x}^H$ has more than one edge.

\begin{thm} \label{2planes}  Take $\xi \in \cM(U_x)$.  Then
$\xi \in \im(\phi)$ if and only if $\xi^H \in \im(\phi^H)$ for all $H
\in \cH$.
\end{thm}

The ``only if'' direction is trivial.  Note that for two-planes  $H
\notin \cH$, $\phi^H$ is automatically surjective.

Pick a subtorus $\C^* \subset T$ which is contracting near $x$.  Its
Lie algebra is a one-dimensional subspace $\mft_0 \subset \mft$; let
$\mft^\bot_0 \subset \mft^*$ be its annihilator.  Since the action is
contracting, we  have $V_L\not\subset \mft_0^\bot$ for all $L\in U_x$.
Let $A_0 = \Sym(\mft^\bot_0)$; it is a subring  of $A$.  Note that the
set of all possible $\mft_0$ forms an open subset of the set of points
in the projective space $\PP(\mft)$ which are rational with respect to
the lattice of characters.  Thus  $\mft_0$ can be chosen to avoid any
finite collection of vectors.

\begin{lemma} $M_{\bdy x}$ and $\cM(U_x)$ are free $A_0$-modules.
\end{lemma}
\begin{proof} The result for $\cM(U_x)$ is clear from the construction of 
$\cM$.  Since in \S4.1 we showed that
\[M_{\bdy x} \cong
IH^*_T(N_0) \cong IH^*_{T/\C^*}(N_0/\C^*),\] we can apply the first
part of Theorem \ref{GKMloc}.
\end{proof}

Now take $\xi \in \cM(U_x)$.  Define an ideal $I(\xi)$ in $A_0$ by
\[I(\xi) = \{a \in A_0 \mid a\xi\in \im(\phi)\}.\]
The previous lemma plus the injectivity of $M_{\bdy x} \to \cM(U_x)$
implies the following.

\begin{prop} (Chang and Skjelbred \cite{CS}) The ideal $I(\xi)$ is principal.
\end{prop}

Take a vector space $H \in \cH\cup \{\{0\}\}$.   We say a vector $v\in
\mft_0^\bot$ is {\em $H$-good} if $v\notin H$, and if, in the case $H
= \{0\}$, $v$ is in some plane $J\in \cH$.

\begin{lemma} \label{a lemma} If $\xi^H \in \im(\phi^H)$, then there is a 
nonzero $p \in I(\xi)$ which is a product of $H$-good linear factors.
\end{lemma}
Note that the condition $\xi^{\{0\}} \in \im(\phi^{\{0\}})$ is vacuous.

Theorem \ref{2planes} immediately follows from the lemma: a generator
of $I(\xi)$ must be a product of linear factors, but if $\xi \in
\im(\phi^H)$ for all $H \in \cH$, none of the possible factors can
actually occur, and so $\xi \in \im(\phi)$.

Before proving the lemma, we need the following easy consequence of
the projectivity of $X$.  We say a moment graph $\Gamma$ is {\em
flexible} at $x$ if for any $H \subset \mft^*$ and any $y\in
\Gamma^H$, $y\ne x$, there is a degree two section $\zeta\in
\cA(\Gamma^H)_2$ so that $\zeta_x = 0$, $\zeta_y\ne 0$, and $\zeta_z
\in H$ for all vertices $z\in \Gamma^H$.

\begin{prop}\label{flexible} The moment graph $\Gamma$ of a projective variety  
is flexible at all its vertices.
\end{prop} 

\begin{proof} The moment map gives an embedding $\mu$ of the vertices of $\Gamma$ into $\mft^*$
so that if $z$ and $w$ are joined by an edge $L$,  then $\mu(z) -
\mu(w)$ is a nonzero vector in $V_L$.  If we choose a linear
projection $p\colon \mft^* \to H$ which does not kill  $\mu(y) -
\mu(x)$, then letting $\zeta_z = p(\mu(z) - \mu(x))$ provides the
required section.
\end{proof}

\begin{proof}[Proof of Lemma \ref{a lemma}]
Let $\wt\Gamma$ be the union of 
$U_x \cup \Gamma^H_{>x}$ with the set of upper vertices of 
edges in $U_x$, and let $\wt{\xi}$ 
be any extension of $\xi$ to $\wt\Gamma$.
We will construct an element $\wt{p} \in \cA(\wt{\Gamma})$ so that
\begin{enumerate}
\item
$\wt{p}|_{U_x}$ comes from an element $p\in A_0$ which is a product 
of $H$-good linear factors, and 
\item for any vertex $y\in \wt{\Gamma}$ and any 
adjacent edge $L \notin \wt{\Gamma}$, $\wt{p}_y \in A\cdot V_L$.
\end{enumerate} 
If we can do this, $\wt{p}\wt{\xi}$ can be extended to $\Gamma_{>x}$ by 
placing a $0$ on all vertices and edges outside of $\wt{\Gamma}$, and so 
$p \in I(\xi)$, as claimed.

Assume that we have chosen $\mft_0$ so that the lines $H \cap \mft_0^\bot$ for 
$H\in \cH$ are all distinct.
 Pick a vertex $y\in \wt{\Gamma}$ and an 
adjacent edge $L \notin \wt{\Gamma}$.  
We will construct a degree two section $a \in \cA(\wt{\Gamma})_2$
satisfying property (1) above and for which $a_y \in V_L$. 
The section $\wt{p}$ we want is the product of  these sections over all choices of
$y$ and $L$.

If $y \notin \Gamma^H_{>x}$, then $y$ is the upper
vertex of an edge $L' \in U_x$.  Since $L' \not\subset \mft_0^\bot$, 
there are nonzero vectors $v \in V_L$, $v' \in \mft_0^\bot$ 
with $v - v' \in V_{L'}$.  The section which is $v$ on $y$ and
$v'$ everywhere else does the trick.  Note that $V_L+V_{L'} \in \cH$, so 
$v'$ lies in a plane in $\cH$.

Now suppose $y\in \Gamma^H_{>x}$, so $H \ne \{0\}$. 
Let $\zeta \in \cA(\Gamma^H)_2$ be the section guaranteed by Proposition
\ref{flexible}.  It extends to a section $\zeta \in \cA(\wt\Gamma)_2$ which is 
zero on all of $U_x$. 
  We can assume that $\mft_0$ 
has been chosen so $\zeta_y \notin\mft_0^\bot$.  Thus we can find $v \in V_L$
so that $v' = v - \zeta_y \in \mft_0^\bot$, and putting $a = v' + \zeta$ gives the 
required section.
\end{proof}

\end{document}